\documentclass[12pt]{amsart}
\input{xypic}
\usepackage{amsmath,amssymb,latexsym,enumerate}

\newtheorem{theorem}{Theorem}[section]
\newtheorem{theoremint}{Theorem}
\newtheorem{propint}[theoremint]{Proposition}
\newtheorem{lemint}[theoremint]{Lemma}
\newtheorem{corint}[theoremint]{Corollary}
\newtheorem{proposition}[theorem]{Proposition}
\newtheorem{propdef}[theorem]{Proposition-Definition}
\newtheorem{corollary}[theorem]{Corollary}
\newtheorem{lemma}[theorem]{Lemma}
\theoremstyle{definition}
\newtheorem{question}{Question}
\newtheorem{definition}[theorem]{Definition}

\newtheorem{remark}[theorem]{Remark}
\newtheorem*{ack}{Acknowledgment}

\newcommand{\Z}{{\mathbf Z}}
\newcommand{\N}{{\mathbf N}}
\newcommand{\R}{{\mathbf R}}
\newcommand{\loto}{\longrightarrow}
\newcommand{\precc}{\preccurlyeq}
\newcommand{\Out}{\textnormal{Out}}
\newcommand{\Aut}{\textnormal{Aut}}

\begin{document}
\title[Isolated points in the space of groups]{On the isolated points in the space of groups}

\author[Cornulier]{Yves de Cornulier}
\address{Yves de Cornulier \\ \'Ecole Polytechnique F\'ed\'erale de Lausanne (EPFL)\\
Institut de G\'eom\'etrie, Alg\`ebre et Topologie (IGAT)\\
CH-1015 Lausanne, Switzerland} \email{yves.cornulier@ens.fr}

\author[Guyot]{Luc Guyot}
\address{Luc Guyot \\ Universit\'e de Gen\`eve \\
Section de math\'ematiques CP 64 \\ 1211 Gen\`eve 4, Swizerland}
\email{Luc.Guyot@math.unige.ch}

\author[Pitsch]{Wolfgang Pitsch}
\address{Wolfgang Pitsch \\ Universitat Aut\`onoma de Barcelona \\
Departament de Matem\`atiques\
\
E-08193 Bellaterra, Spain}
\email{pitsch@mat.uab.es}

\thanks{The second author is supported by the Swiss National Science
  Foundation, No.~PP002-68627. The third author is supported by the MEC
  grant MTM2004-06686 and by the program
Ram\'on y Cajal, MEC, Spain}
%    General info

%20E05 Free nonabelian groups
%20E07 Subgroup theorems; subgroup growth
%20E15 Chains and lattices of subgroups, subnormal subgroups
%20E22 Extensions, wreath products, and other compositions
%20E26 Residual properties and generalizations
%20E34 General structure theorems
%20F05 Generators, relations, and presentations
%20F10 Word problems, other decision problems, connections with logic and automata
%20F65 Geometric group theory
%57M07 Topological methods in group theory

\date{Mar. 8, 2006}

%%%%%%%%%%%%%%%%%%%%%%%%%%%%%%%%%%%%%%%%%%%%%%%%%%%%%%%%%%%%%%%%%
%                         ABSTRACT                              %
%%%%%%%%%%%%%%%%%%%%%%%%%%%%%%%%%%%%%%%%%%%%%%%%%%%%%%%%%%%%%%%%%

\begin{abstract}We investigate the isolated points in the
space of finitely generated groups. We give a workable
characterization of isolated groups and study their hereditary
properties. Various examples of groups are shown to yield isolated
groups. We also discuss a connection between isolated groups and
solvability of the word problem.
\end{abstract}

\maketitle

%%%%%%%%%%%%%%%%%%%%%%%%%%%%%%%%%%%%%%%%%%%%%%%%%%%%%%%%%%%%%%%%%
%                       INTRODUCTION                            %
%%%%%%%%%%%%%%%%%%%%%%%%%%%%%%%%%%%%%%%%%%%%%%%%%%%%%%%%%%%%%%%%%

%%%%%%%%%%%%%%%%%%%%%%%%%%%%%%%%%%%%%%%%%%%%%%%%%%%%%%%%%%%%%%%%%
%%%%%%%%%%%%%%%%%%%%%%%%%%%%%%%%%%%%%%%%%%%%%%%%%%%%%%%%%%%%%%%%%
\section*{Introduction}
%%%%%%%%%%%%%%%%%%%%%%%%%%%%%%%%%%%%%%%%%%%%%%%%%%%%%%%%%%%%%%%%%
%%%%%%%%%%%%%%%%%%%%%%%%%%%%%%%%%%%%%%%%%%%%%%%%%%%%%%%%%%%%%%%%%

At the end of his celebrated paper ``Polynomial growth and
expanding maps'' \cite{MR0623534}, Gromov sketched what could be a
topology on a set of groups. His ideas led to the construction by
Grigorchuk of the space of marked groups \cite{MR0764305}, where
points are finitely generated groups with $m$ marked generators.
This ``space of marked groups of rank $m$" $\mathcal{G}_m$ is a
totally discontinuous compact metrizable space.

One of the main interests of $\mathcal{G}_m$ is to find properties
of groups that are reflected in its topology. Various elementary
observations in these directions are made in \cite{MR2151593}: for
instance, the class of nilpotent finitely generated groups is
open, while the class of solvable finitely generated groups is
not; the class of finitely generated orderable groups is closed,
etc. Deeper results can be found about the closure of free groups
(see \cite{MR2151593} and the references therein) and the closure
of hyperbolic groups \cite{MR1760424}. The study of the
neighbourhood of the first Grigorchuk group has also proved to be
fruitful \cite{MR0764305} in the context of growth of finitely
generated groups.

An example of an open question about $\mathcal{G}_m$ ($m\ge 2$) is
the following: does there exist a surjective continuous invariant:
$\mathcal{G}_m\to [0,1]$? On the other hand, it is known that
there exists no real-valued injective measurable invariant
\cite{MR1760424}.

The aim of this paper is the study of the isolated points in
$\mathcal{G}_m$, which we call \textit{isolated groups}. It turns
out that these groups have already occurred in a few papers
\cite{MR0414671,MR0647191}, without the topological point of view.
They are introduced by B.H.~Neumann \cite{MR0414671} as ``groups
with finite absolute presentation". It follows from a result of
Simmons \cite{MR0342614} that they have solvable word problem, see
the discussion in \S\ref{sec word}. However the only examples
quoted in the literature are finite groups and finitely presented
simple groups; we provide here examples showing that the class of
isolated groups is considerably larger.

Let us now describe the paper. In Section \ref{sec spacegroups} we
construct the space of finitely generated groups.

Here is an elementary but useful result about this topology:
\begin{lemint}\label{lem local_indep_intro}
Consider two marked groups $G_1\in\mathcal{G}_{m_1}$,
$G_2\in\mathcal{G}_{m_2}$. Suppose that they are isomorphic. Then
there are clopen (= closed open) neighbourhoods $V_i$, $i=1,2$ of
$G_i$ in $\mathcal{G}_{m_i}$ and a homeomorphism $\varphi:V_1\to
V_2$ mapping $G_1$ to $G_2$ and preserving isomorphism classes,
i.e. such that, for every $H\in V_1$, $\varphi(H)$ is isomorphic
to $H$ (as abstract groups).
\end{lemint}

This allows us to speak about the ``space of finitely generated
groups\footnote{This must \textit{not} be viewed as the set of
isomorphism classes of finitely generated groups, even locally.
Indeed, an infinite family of isomorphic groups may accumulate at
the neighbourhood of a given point (e.g. groups isomorphic to $\Z$
at the neighbourhood of $\Z^2$ \cite[Example~2.4(c)]{MR2151593}),
and we certainly do not identify them. What we call ``space of
finitely generated groups" might be viewed as the disjoint union
of all $\mathcal{G}_n$, but this definition is somewhat arbitrary;
and is not needed since we only consider local properties.}"
rather than ``space of marked groups" whenever we study local
topological properties.
In particular, to be isolated is an algebraic property of the
group, i.e. independent of the marking. This bears out the
terminology of ``isolated groups".

In Section \ref{sect isolated}, we proceed with the
characterization of isolated groups. In a group $G$, we call a
subset $F$ a \textit{discriminating subset} if every non-trivial
normal group of $G$ contains an element of $F$, and we call $G$
\textit{finitely discriminable} if it has a \textit{finite}
discriminating subset. Finitely discriminable groups are
introduced as ``semi-monolithic groups" in \cite{MR0647191}.
Here is an algebraic characterization of isolated groups (compare
with \cite[Proposition~2(a)]{MR0647191}).

%%%%%%%%%%%%%%%%%%%%%%%%%%%%%%%%%%%%%%%%%%%%%%%%%%%%%%%%%%%%%%%%%%
\begin{propint}\label{thm intro}
%%%%%%%%%%%%%%%%%%%%%%%%%%%%%%%%%%%%%%%%%%%%%%%%%%%%%%%%%%%%%%%%%%
A group $G$ is isolated if and only if the two following
properties are satisfied
\begin{enumerate}[(i)]
\item $G$ finitely presented; \item $G$ is finitely discriminable.
\end{enumerate}
\end{propint}
%%%%%%%%%%%%%%%%%%%%%%%%%%%%%%%%%%%%%%%%%%%%%%%%%%%%%%%%%%%%%%%%%%

Since the class of finitely presented groups is well understood in
many respects, we are often led to study finite discriminability.

%%%%%%%%%%%%%%%%%%%%%%%%%%%%%%%%%%%%%%%%%%%%%%%%%%%%%%%%%%%%%%%%%%
\begin{propint}
%%%%%%%%%%%%%%%%%%%%%%%%%%%%%%%%%%%%%%%%%%%%%%%%%%%%%%%%%%%%%%%%%%
The subspace of finitely discriminable groups is dense in the
space of finitely generated groups.
\end{propint}
%%%%%%%%%%%%%%%%%%%%%%%%%%%%%%%%%%%%%%%%%%%%%%%%%%%%%%%%%%%%%%%%%%
In other words, the subspace of finitely discriminable groups is
dense in $\mathcal{G}_m$ for all $m \geq 1$.

In Section \ref{sec word}, we discuss the connection with the word
problem. We call a finitely generated group $G$ recursively
discriminable if there exists in $G$ a recursively enumerable
discriminating subset (see Section \ref{sec word} for a more
precise definition if necessary).
Proposition \ref{thm intro} has the following analog
\cite[Theorem~B]{MR0342614}.
%%%%%%%%%%%%%%%%%%%%%%%%%%%%%%%%%%%%%%%%%%%%%%%%%%%%%%%%%%%%%%%%%%
\begin{theoremint}
%%%%%%%%%%%%%%%%%%%%%%%%%%%%%%%%%%%%%%%%%%%%%%%%%%%%%%%%%%%%%%%%%%
A finitely generated group has solvable word problem if and only
if it is both recursively presentable and recursively
discriminable.\end{theoremint}
%%%%%%%%%%%%%%%%%%%%%%%%%%%%%%%%%%%%%%%%%%%%%%%%%%%%%%%%%%%%%%%%%%

This is a conceptual generalization of a well-known theorem by
Kuznetsov stating that a recursively presentable simple group has
solvable word problem.

\begin{corint}
An isolated group has solvable word problem.\qed
\end{corint}

The existence of certain pathological examples of finitely
presented groups due to Miller III \cite{MR0629262} has the
following consequence.

\begin{propint}
The class of groups with solvable word problem is not dense in the
space of finitely generated groups. In particular, the class of
isolated groups is not dense.
\end{propint}

%%%%%%%%%%%%%%%%%%%%%%%%%%%%%%%%%%%%%%%%%%%%%%%%%%%%%%%%%%%%%%%%%%%%%%%%%%%%%
In Section \ref{sec hereditary} we %use the results of Section
%\ref{sect isolated} to
explore the hereditary properties of
finitely discriminable groups, and thus isolated groups. For instance we prove

%%%%%%%%%%%%%%%%%%%%%%%%%%%%%%%%%%%%%%%%%%%%%%%%%%%%%%%%%%%%%%%%%%
\begin{theoremint}\label{thm intro_extensions}
%%%%%%%%%%%%%%%%%%%%%%%%%%%%%%%%%%%%%%%%%%%%%%%%%%%%%%%%%%%%%%%%%%
The class of finitely discriminable (resp. isolated) groups is
stable under:
\begin{enumerate}
\item[1)] extensions of groups; \item[2)] taking overgroups of
finite index, i.e. if $H$ is finitely discriminable (resp.
isolated) and of finite index in $G$ then so is $G$.
\end{enumerate}
\end{theoremint}
%%%%%%%%%%%%%%%%%%%%%%%%%%%%%%%%%%%%%%%%%%%%%%%%%%%%%%%%%%%%%%%%%%

The proof of Theorem \ref{thm intro_extensions} is less immediate
than one might expect. For instance it involves the classification
of finitely discriminable abelian groups (Lemma \ref{lem centerofisolated}). We state a much more
general result in Theorem \ref{thm stabextension} (see also
\S\ref{subs:wreath} for the case of wreath products).

Finally in Section \ref{sect exaisolated} we provide examples of
isolated groups. First note that the most common infinite finitely
generated groups are not finitely discriminable. For instance, if
$G$ is an infinite residually finite group, then for every finite
subset $F\subset G-\{1\}$ there exists a normal subgroup of finite
index $N$ satisfying $N\cap F=\emptyset$. This prevents $G$ from
being finitely discriminable. On the other hand, Champetier
\cite{MR1760424} has proved that the closure in $\mathcal{G}_n$
($n\ge 2$) of the set of non-elementary hyperbolic groups, is a
Cantor set and therefore contains no isolated point. We leave for
record

\begin{propint}
Infinite residually finite groups and infinite hyperbolic groups
are not finitely discriminable.\qed
\end{propint}

On the other hand, the simplest examples of isolated groups are
finite groups. There are also finitely presented simple groups.
But the class of isolated groups is considerably larger, in view
of the following result, proved in \S\ref{subs:quotient_isolated}:

%%%%%%%%%%%%%%%%%%%%%%%%%%%%%%%%%%%%%%%%%%%%%%%%%%%%%%%%%%%%%%%%%%
\begin{theoremint}
Every finitely generated group is a quotient of an isolated
group.\end{theoremint}
%%%%%%%%%%%%%%%%%%%%%%%%%%%%%%%%%%%%%%%%%%%%%%%%%%%%%%%%%%%%%%%%%%

This shows in particular that the lattice of normal subgroups of
an isolated group can be arbitrarily complicated; for instance, in
general it does not satisfy the descending/ascending chain
condition.

\begin{propint}
There exists an isolated group that is 3-solvable and
non-Hopfian\footnote{A group $G$ is non-Hopfian if there exists an
epimorphism $G\to G$ with non-trivial kernel.}.\label{thm
isol_solv_intro}
\end{propint}

This is in a certain sense optimal, since it is known that
finitely generated groups that are either nilpotent or metabelian
(2-solvable) are residually finite \cite{MR0110750}, and thus
cannot be isolated unless they are finite. The example we provide
to prove Proposition \ref{thm isol_solv_intro} (see
\S\ref{subs:abels}) is a group that had been introduced by Abels
\cite{MR564423} as the first example of a non-residually finite
(actually non-Hopfian) finitely presented solvable group. A
variation on this example provides (recall that a countable group
$G$ has Kazhdan's Property~T if every isometric action of $G$ on a
Hilbert space has a fixed point):

\begin{propint}
There exists an infinite isolated group with Kazhdan's Property~T.
\end{propint}

We provide some other examples. One of them (see
\S\ref{subs:houghton}) is Houghton's group, which is an
extension of the group of finitely supported permutations of a
countable set, by $\Z^2$. In particular, this group is elementary
amenable but non virtually solvable.

Another one (see \S\ref{subs:grigorchuk_fp}) is a group exhibited
by Grigorchuk \cite{MR1616436}, which is the first known finitely
presented amenable group that is not elementary amenable. This is
an ascending HNN extension of the famous ``first Grigorchuk
group", which has intermediate growth and is torsion; the latter
is certainly not isolated since it is infinitely presented and is
not finitely discriminable since it is residually finite. The fact
that this group is isolated contradicts a conjecture by Stepin in
\cite[\S 1]{MR1616436}, stating that every amenable finitely
generated group can be approximated by elementary amenable ones.

Finally (see \S\ref{subs:deligne}), some lattices in non-linear
simple Lie groups provide examples of isolated groups that are
extensions with infinite residually finite quotient and finite
central kernel.

\medskip

Throughout this article we use the following notation. If $x$ and
$y$ are elements in a group $G$ then
\[
[x,y] = xyx^{-1}y^{-1}, x^y = y^{-1}xy.
\]
Similarly if $N$ and $K$ are subgroups of $G$, then $[N,K]$ stands
for the subgroup generated by $\{ [n,k] \ | \ n \in N, k \in K\}$.
Finally in any group $G$, we denote by $Z(G)$ the center, and more generally by $C_G(X)$ the centralizer of a subset $X\subset G$.

\begin{ack}We thank Avinoam Mann for pointing out some useful
references.
\end{ack}

%%%%%%%%%%%%%%%%%%%%%%%%%%%%%%%%%%%%%%%%%%%%%%%%%%%%%%%%%%%%%%%%%%%%%%%
%%%%%%%%%%%%%%%%%%%%%%%%%%%%%%%%%%%%%%%%%%%%%%%%%%%%%%%%%%%%%%%%%%%%%%%
\section{The space of finitely generated groups}\label{sec spacegroups}
%%%%%%%%%%%%%%%%%%%%%%%%%%%%%%%%%%%%%%%%%%%%%%%%%%%%%%%%%%%%%%%%%%%%%%%
%%%%%%%%%%%%%%%%%%%%%%%%%%%%%%%%%%%%%%%%%%%%%%%%%%%%%%%%%%%%%%%%%%%%%%%

Let \(G\) be a group. We denote by \(\mathcal{P}(G)\) the set of
subsets of \(G\) and by \(\mathcal{G}(G)\) the set of normal
subgroups of \(G\). We endow $\mathcal{P}(G)$ with the product
topology through the natural bijection with $\{0,1\}^G$. Hence
\(\mathcal{P}(G)\) is a compact and totally discontinuous space.

%By definition, a sequence \(A_n\)  converges to \(A\) in
%\(\mathcal{P}(G)\) if and only if for any finite subset \(B\) of
%\(G\) there is an integer \(N\) such that \(B \cap A_n= B \cap A\)
%for all \(n \ge N\).

Limits in $\mathcal{P}(G)$ have the following simple description,
whose proof is straightforward and omitted.

%%%%%%%%%%%%%%%%%%%%%%%%%%%%%%%%%%%%%%%%%%%%%%%%%%%%%%%%%%%%%%%%%
\begin{lemma}\label{lem descrlimits}
%%%%%%%%%%%%%%%%%%%%%%%%%%%%%%%%%%%%%%%%%%%%%%%%%%%%%%%%%%%%%%%%%
The net \(( A_i) \) converges to  $A$ in \(\mathcal{P}(G)\) if and
only if \( A = \liminf A_i = \limsup A_i \), where \( \liminf A_i
= \bigcup_{i} \bigcap_{j \ge i}A_j\) and \( \limsup A_i =
\bigcap_{i} \bigcup_{j \ge i} A_i\).\qed
\end{lemma}
%%%%%%%%%%%%%%%%%%%%%%%%%%%%%%%%%%%%%%%%%%%%%%%%%%%%%%%%%%%%%%%%%

Since, for every net $(N_i)$ of normal subgroups, $\liminf N_i$ is
also a normal subgroup, the following proposition follows
immediately from the lemma.

%%%%%%%%%%%%%%%%%%%%%%%%%%%%%%%%%%%%%%%%%%%%%%%%%%%%%%%%%%%%%%%%%
\begin{proposition} \label{prop normalclosed}
%%%%%%%%%%%%%%%%%%%%%%%%%%%%%%%%%%%%%%%%%%%%%%%%%%%%%%%%%%%%%%%%%
The subset \(\mathcal{G}(G)\) is closed in \(\mathcal{P}(G)\).\qed
\end{proposition}
%%%%%%%%%%%%%%%%%%%%%%%%%%%%%%%%%%%%%%%%%%%%%%%%%%%%%%%%%%%%%%%%%

The space $\mathcal{G}(G)$ can be identified with the space of
quotients of $G$, which we also call, by abuse of notation,
$\mathcal{G}(G)$ (in the sequel it will always be clear when we
consider an element of $\mathcal{G}$ as a normal subgroup or as a
quotient of $G$). It is endowed with a natural order: $H_1\precc
H_2$ if the corresponding normal subgroups $N_1,N_2$ satisfy
$N_1\supset N_2$.

The topology of $\mathcal{G}(G)$ has the following basis:
$$(\Omega_{r_1,\dots,r_k,s_1,\dots,s_\ell}),\quad
r_1,\dots,r_k,s_1,\dots,s_l\in G,$$ where
$\Omega_{r_1,\dots,r_k,s_1,\dots,s_\ell}$ is the set of quotients
of $G$ in which each $r_i=1$ and each $s_j\neq 1$. These are open
and closed subsets.

If $F_m$ is a free group of rank $m$ with a given freely
generating family, then $\mathcal{G}(F_m)$ is usually called the
\emph{space of marked groups on $m$ generators} and we denote it
by $\mathcal{G}_m$. An element in $\mathcal{G}_m$ can be viewed as
a pair $(G,\mathcal{T})$ where $G$ is an $m$-generated group and
$\mathcal{T}$ is a generating $m$-tuple.

If $G,H$ are any groups, every homomorphism $f:G\to H$ induces a
continuous map $f^*:\mathcal{G}(H)\to\mathcal{G}(G)$, which is
injective if $f$ is surjective. The main features of the spaces
$\mathcal{G}(G)$ is summarized in the following Lemma, which is
essentially known (see \cite[Lemme~2.2 and
Proposition~3.1]{MR1760424}).

%%%%%%%%%%%%%%%%%%%%%%%%%%%%%%%%%%%%%%%%%%%%%%%%%%%%%%%%%%%%%%%%%%
\begin{lemma}~
%%%%%%%%%%%%%%%%%%%%%%%%%%%%%%%%%%%%%%%%%%%%%%%%%%%%%%%%%%%%%%%%%%
\begin{itemize}
\item[(1)] Let $G$ be a group and $H$ a quotient of $G$; denote by
$p$ the quotient map $G\to H$. Then the embedding
$p^*:\mathcal{G}(H)\to\mathcal{G}(G)$ is a closed homeomorphism
onto its image, which we identify with $\mathcal{G}(H)$. Moreover,
the following are equivalent.
    \begin{itemize}
        \item[(i)] $\mathcal{G}(H)$ is open in
        $\mathcal{G}(G)$.
        \item[(ii)] $H$ is contained in the interior of
        $\mathcal{G}(H)$ in $\mathcal{G}(G)$ (in other words: $H$ has a
        neighbourhood in $\mathcal{G}(G)$ consisting of quotients of itself).
        \item[(iii)] $\textnormal{Ker}(p)$ is finitely generated
        as a normal subgroup of $G$.
    \end{itemize}
\item[(2)] If, in addition, $G$ is a finitely presented group,
then these are also equivalent to
    \begin{itemize}
    \item[(iv)] $H$ is finitely presented.
    \end{itemize}
\item[(3)] Let $G_1$, $G_2$ be finitely presented groups, and
consider quotients $H_i\in\mathcal{G}(G_i)$, $i=1,2$. Suppose that
$H_1$ and $H_2$ are isomorphic groups. Then there exist finitely
presented intermediate quotients $H_i\precc K_i\precc G_i$,
$i=1,2$, and an isomorphism $\phi:K_1\to K_2$, such that $\phi^*$
maps the point $H_2$ to $H_1$.
    \end{itemize}
\label{lem:nei}
\end{lemma}

Note that Lemma \ref{lem local_indep_intro} is an immediate
consequence of Lemma \ref{lem:nei}.

%%%%%%%%%%%%%%%%%%%%%%%%%%%%%%%%%%%%%%%%%%%%%%%%%%%%%%%%%%%%%%%%%%%
\begin{proof}(1) is straightforward and left to the reader. (2)
follows from (1) and the fact that if $G$ is a finitely presented
group and $N$ a normal subgroup, then $G/N$ is finitely presented
if and only if $N$ is finitely generated as a normal subgroup
\cite[Lemma~11.84]{MR1307623}.

To prove (3), fix an isomorphism $\alpha:H_1\to H_2$ and identify
$H_1$ and $H_2$ to a single group $H$ through $\phi$. Take a
finitely generated subgroup $F_0$ of the fibre product
$G_1\times_H G_2$ mapping onto both $G_1$ and $G_2$, and take a
free group $F$ of finite rank mapping onto $F_0$. Then the
following diagram commutes.
\[
\xymatrix{ F \ar@{->}[r] \ar@{->}[d]& G_2\ar@{->}[d] \\
G_1\ar@{->}[r] & H }
\]%

For $i=1,2$, let $N_i$ denote the kernel of $F\to G_i$. Then
$F/N_1N_2$ is a finitely presented quotient of both $G_1$ and
$G_2$, having $H$ as a quotient. If we do not longer identify
$H_1$ and $H_2$, then $F/N_1N_2$ can be viewed as a quotient $K_i$
of $G_i$, and the obvious isomorphism $\phi$ between $K_1$ and
$K_2$ (induced by the identity of $F$)
induces~$\alpha$.\nolinebreak\end{proof}

It follows that every local topological consideration makes sense
in a ``space of finitely generated groups". Roughly speaking the
later looks like a topological space. Its elements are finitely
generated groups (and therefore do not make up a set). If $G$ is a
finitely generated group, then a neighbourhood of $G$ is given by
$\mathcal{G}(H)$, where $H$ is a finitely presented group endowed
with a given homomorphism onto $G$. For instance, if $\mathcal{C}$
is an isomorphism-closed class of groups (as all the classes of
groups we consider in the paper), then we can discuss whether
$\mathcal{C}$ is open, whether it is dense.

%%%%%%%%%%%%%%%%
%%%%%%%%%%%%%%%%%%%%%%%%%%%%%%%%%%%%%%%%%%%%%%%%%
%%%%%%%%%%%%%%%%%%%%%%%%%%%%%%%%%%%%%%%%%%%%%%%%%%%%%%%%%%%%%%%%%
\section{Isolated groups}\label{sect isolated}
%%%%%%%%%%%%%%%%%%%%%%%%%%%%%%%%%%%%%%%%%%%%%%%%%%%%%%%%%%%%%%%%%
%%%%%%%%%%%%%%%%%%%%%%%%%%%%%%%%%%%%%%%%%%%%%%%%%%%%%%%%%%%%%%%%%
%
%%%%%%%%%%%%%%%%%%%%%%%%%%%%%%%%%%%%%%%%%%%%%%%%%%%%%%%%%%%%%%%%%%%%%%%%%%%%%%%thm>prop
\begin{propdef}\label{thm finimin}
%%%%%%%%%%%%%%%%%%%%%%%%%%%%%%%%%%%%%%%%%%%%%%%%%%%%%%%%%%%%%%%%%%%%%%%%%%%%%%%
Let $G$ be a group. We say that $G$ is {\em finitely
discriminable} if it satisfies the following equivalent
conditions.
\begin{enumerate}
    \item[i)] The trivial normal subgroup $\{1\}$ is isolated in $\mathcal{G}(G)$.
    \item[ii)] The group $G$ has finitely many minimal normal
    subgroups and any non-trivial normal subgroup contains a minimal one.
    \item[iii)] There exists a finite {\em discriminating subset} in $G$: this
    is a finite subset $F \subset G-\{1\}$ such
    that any non-trivial normal subgroup of $G$ contains at
    least one element of $F$.
\end{enumerate}
\end{propdef}
%%%%%%%%%%%%%%%%%%%%%%%%%%%%%%%%%%%%%%%%%%%%%%%%%%%%%%%%%%%%%%%%%%%%%%%%%%%%%%%
\begin{proof}[Proof of the equivalences:]
$ii) \Rightarrow iii)$ Define $F$ by taking a non-trivial element
in each minimal normal subgroup.

$iii) \Rightarrow ii)$ Every non-trivial normal subgroup contains
the normal subgroup generated by an element of $F$.

$iii) \Rightarrow i)$ Since $F$ is finite, the set of normal
subgroups with empty intersection with $F$ is open in
$\mathcal{G}(G)$; by assumption this set is reduced to
$\{\{1\}\}$.

$i) \Rightarrow iii)$ We contrapose. For every finite subset
$F\subset G-\{1\}$, consider a non-trivial normal subgroup $N_F$
having empty intersection with $F$. Then $N_F\to\{1\}$ when $F$
becomes large (that is, tends to $G-\{1\}$ in $\mathcal{P}(G)$).
\end{proof}

%%%%%%%%%%%%%%%%%%%%%%%%%%%%%%%%%%%%%%%%%%%%%%%%%%%%%%%%%%%%%%%%%%%
\begin{proposition}\label{propcar}
%%%%%%%%%%%%%%%%%%%%%%%%%%%%%%%%%%%%%%%%%%%%%%%%%%%%%%%%%%%%%%%%%%%
A group $G$ is isolated if and only if it is both finitely
presentable and finitely discriminable.
\end{proposition}
%%%%%%%%%%%%%%%%%%%%%%%%%%%%%%%%%%%%%%%%%%%%%%%%%%%%%%%%%%%%%%%%%%%%
\begin{proof}
An isolated group is finitely presented: this follows from the
implication (ii)$\Rightarrow$(iv) of Lemma \ref{lem:nei}.

Suppose that $G$ is finitely generated but not finitely
discriminable. Then $G$ is not isolated in $\mathcal{G}(G)$
(viewed as the set of quotients of $G$), hence is not isolated in
the space of finitely generated groups.

Conversely suppose that $G$ satisfies the two conditions. Since
$G$ is finitely presented, by the implication (iv)$\Rightarrow$(i)
of Lemma \ref{lem:nei}, $\mathcal{G}(G)$ is a neighbourhood of $G$
in the space of finitely generated groups. Since $G$ is finitely
discriminable, it is isolated in $\mathcal{G}(G)$. Hence $G$ is
isolated.\nolinebreak\end{proof}

Proposition \ref{propcar} allows us to split the study of isolated
groups into the study of finitely discriminable groups and
finitely presented groups. These studies are in many respects
independent, and it is sometimes useful to drop the finite
presentability assumption when we have to find examples. However,
these properties also have striking similarities, for instance:

%%%%%%%%%%%%%%%%%%%%%%%%%%%%%%%%%%%%%%%%%%%%%%%%%%%%%%%%%%%%%%%%%%%%%%%%%%%%%%%%thm->prop
\begin{proposition}\label{thm fddense}
The classes of finitely discriminable and finitely presentable
groups are both dense in the space of finitely generated groups.
\end{proposition}
\begin{proof}
The case of finitely presentable groups is an observation by
Champetier \cite[Lemme~2.2]{MR1760424} (it suffices to approximate
every finitely generated group by truncated presentations).

Let us deal with finite discriminability. If $G$ is finite then it
is finitely discriminable. Otherwise, for every finite subset $F$
of $G-\{1\}$, consider a maximal normal subgroup $N_F$ among those
with empty intersection with $F$. Then $G/N_F$ is finitely
discriminated by the image of $F$, while the net $(G/N_F)$
converges to $G$ when $F$ becomes large.
\end{proof}

\begin{remark}
The proof above shows that more generally, every group $G$ is
approximable by finitely discriminable quotients in
$\mathcal{G}(G)$.
\end{remark}

In contrast, we show in the next paragraph that being isolated is
not a dense property.

%%%%%%%%%%%%%%%%%%%%%%%%%%%%%%%%%%%%%%%%%%%%%%%%%%%%%%%%%%%%%%%%%
%%%%%%%%%%%%%%%%%%%%%%%%%%%%%%%%%%%%%%%%%%%%%%%%%%%%%%%%%%%%%%%%%
\section{Isolated groups and the word problem}\label{sec word}
%%%%%%%%%%%%%%%%%%%%%%%%%%%%%%%%%%%%%%%%%%%%%%%%%%%%%%%%%%%%%%%%%
%%%%%%%%%%%%%%%%%%%%%%%%%%%%%%%%%%%%%%%%%%%%%%%%%%%%%%%%%%%%%%%%%

Roughly speaking, a sequence of words in a free group $F$ of
finite rank is $\textit{recursive}$ if it can be computed by
a finite algorithm; more precisely by an ideal computer, namely a
Turing machine. See Rotman's book \cite[Chap.~12]{MR1307623} for a
precise definition. A subset $X\subset F$ is \textit{recursively
enumerable} if it is the image of a recursive sequence.

\begin{propdef}
Let $G$ be a finitely generated group. We call a sequence $(g_n)$
in $G$ \textit{recursive} if it satisfies one of the two
equivalent properties:

\begin{itemize}
\item[(1)] There exists a free group of finite rank $F$, an
epimorphism $p:F\to G$, and a recursive sequence $(h_n)$ in $F$
such that $p(h_n)=g_n$ for all $n$.
 \item[(2)] For every free group of finite rank $F$ and every epimorphism
$p:F\to G$, there exists a recursive sequence $(h_n)$ in $F$ such
that $p(h_n)=g_n$ for all $n$.
\end{itemize}
\end{propdef}
\begin{proof}We have to justify that the two conditions are
equivalent. Note that (2) is a priori stronger. But if (1) is
satisfied, then, using Tietze transformations to pass from a
generating subset of $G$ to another, we obtain that (2) is
satisfied.
\end{proof}

\begin{definition}
A finitely generated group $G$ is \emph{recursively discriminable}
if there exists a recursively enumerable discriminating subset:
there exists a recursive sequence $(g_n)$ in $G-\{1\}$ such that
every normal subgroup $N\neq 1$ of $G$ contains some~$g_n$.
\end{definition}

\begin{remark}
A related notion, namely that of terminal groups, is introduced by
A.~Mann in \cite[Definition~2]{MR0647191}. A finitely generated
group $G=F/N$, with $F$ free of finite rank, is \textit{terminal}
if $F-N$ is recursively enumerable. Clearly this implies that $G$
is recursively discriminable, but the converse is false: indeed,
there are only countably many terminal groups, while there are
$2^{\aleph_0}$ non-isomorphic finitely generated simple groups
\cite[Chap.~IV, Theorem~3.5]{MR0577064}. Nevertheless, there exist
terminal groups with unsolvable word problem
\cite[Proposition~1]{MR0647191}.
\end{remark}

Let $G$ be a finitely generated group, and write $G=F/N$ with $F$
a free group of finite rank. Recall that $G$ is recursively
presentable if and only if $N$ is recursively enumerable, and that
$G$ has solvable word problem if and only if both $N$ and $F-N$
are recursively enumerable; that is, $N$ is recursive.

The following theorem was originally proved by Simmons
\cite{MR0342614}. We offer here a much more concise proof of this
result. It can be viewed as a conceptual generalization of a
well-known theorem of Kuznet\-sov \cite[Chap.~IV,
Theorem~3.6]{MR0577064}, which states that a recursively
presentable simple group has solvable word problem.

\begin{theorem}
Let $G=F/N$ be a finitely generated group as above. Then $G$ has
solvable word problem if and only if it is both recursively
presentable and recursively discriminable. \label{thm word_pb}
\end{theorem}
\begin{proof}
The conditions are clearly necessary. Conversely, suppose that
they are satisfied, and let us show that $F/N$ has solvable word
problem. Consider a recursive discriminating sequence $(g_n)$ in
$F-\nolinebreak N$.

Let $x$ belong to $F$, and let $N_x$ be the normal subgroup
generated by $N$ and $x$; it is recursively enumerable. Set
$W_x=\{y^{-1}g_n\,|\;y\in N_x\;n\in\N\}$. Then $W_x$ is
recursively enumerable. Observe that $1\in W_x$ if and only if
$x\notin N$. Indeed, if $1\in W_x$, then $g_n=y$ for some $y\in
N_x$, and since $g_n\notin N$ this implies that $x\notin N$.
Conversely if $x\notin N$, then $N_x$ projects to a non-trivial
subgroup of $F/N$, and hence contains one of the $g_n$'s. So the
algorithm is the following: enumerate both $N$ and $W_x$: either
$x$ appears in $N$, or $1$ appears in $W_x$ and in this case
$x\notin N$.
\end{proof}

Our initial motivation in proving Theorem \ref{thm word_pb} is the
following corollary.

\begin{corollary}
An isolated group has solvable word problem.\qed
\end{corollary}

There exists an alternative short proof of the corollary using
model Theory (see the proof of the analogous assertion in
\cite{MR0642422}). More precisely, groups with solvable word
problem are characterized by Rips \cite{MR0642422} as isolated
groups for some topology on $\mathcal{G}_n$ stronger than the
topology studied here, making the corollary obvious.

Since the class of isolated groups is the intersection of the
classes of finitely presented and finitely generated finitely
discriminable groups and since these classes are both dense, it is
natural to ask whether the class of isolated groups is itself
dense. This question was the starting point of our study. However
it has a negative answer.

\begin{proposition}
The class of finitely generated groups with solvable word problem
is not dense.\label{prop word_not_dense}
\end{proposition}

\begin{corollary}
The class of isolated groups is not dense.\qed
\end{corollary}

\begin{proof}[Proof of Proposition \ref{prop word_not_dense}]
C.~Miller~III has proved \cite{MR0629262} that there exists a
non-trivial finitely presented group $G$ such that the only
quotient of $G$ having solvable word problem is $\{1\}$ (and
moreover $G$ is SQ-universal, i.e. every countable group embeds in
some quotient of $G$). Thus, using Lemma \ref{lem:nei}, $G$ is not
approximable by isolated groups.\nolinebreak
\end{proof}

This leaves many questions open.

\begin{question}Is every finitely generated group with solvable word problem
a limit of isolated groups?\end{question}

\begin{question}Is every word hyperbolic group a limit of isolated groups?
\end{question}
Note that a word hyperbolic group has solvable word problem. The
following stronger question is open: is every word hyperbolic
group residually finite?

\begin{question}Is every finitely generated solvable group a limit of isolated groups?
\end{question}
Note that there exist finitely presented solvable groups with
unsolvable word problem \cite{MR0631441}; this suggests a negative
answer.

\medskip

Let $G$ be a group. An \textit{equation} (resp.
\textit{inequation}) over $G$ is an expression ``$m=1$" (resp.
$``m\neq 1"$) where $m=m(x_1,\dots,x_n)$ is an element of the free
product of $G$ with the free group over unknowns $x_1,\dots,x_n$.
A solution of an (in)equation in $G$ is a $n$-uple
$(g_1,\dots,g_n)$ of $G$ such that $m(g_1,\dots,g_n)=1$ (resp.
$m(g_1,\dots,g_n)\neq 1$). A system of equations and inequations
over a group $G$ is \textit{coherent} if it has solution in some
overgroup of $G$. A group $\Omega$ is \textit{existentially
closed} if every coherent finite system of equations and
inequations over $\Omega$ has a solution in $\Omega$.

Using free products with amalgamation, one can prove
\cite{MR0040299} that every group $G$ embeds in an existentially
closed group, which can be chosen countable if $G$ is so. The
\textit{skeleton} of a group $G$ is the class of finitely
generated groups embedding in $G$. Skeletons of existentially
closed groups have been extensively studied (see
\cite{MR0960689}). What follows is not new but merely transcribed
in the language of the space of finitely generated groups. Our aim
is to prove that this language is relevant in this context.

The link with isolated groups is given by the following
observation \cite[Lemma~2.4]{MR0414671}: a given isolated group
embeds in every existentially closed group. This is contained in
the following more general result.

\begin{proposition}
If $\Omega$ is an existentially closed group, then its skeleton is
dense in the space of finitely generated groups.\label{prop
skeleton dense}
\end{proposition}
\begin{proof}Fix an existentially closed group $\Omega$.
Let $F$ be a free group of rank $n$, and $g_1,\dots,g_n$ be
generators. Choose elements $m_1,\dots,m_d,\mu_1,\dots,\mu_\delta$
in $F$. Consider the set $\mathcal{S}$ of quotients of $F$ in
which $m_1,\dots,m_d=1$ and $\mu_1,\dots,\mu_\delta\neq 1$.
Subsets of this type make up a basis of open (and closed) subsets
in $\mathcal{G}(F)$. To say that such a subset $\mathcal{S}$ is
non-empty means that the system of equations and inequations
$$\left\{%
\begin{array}{ll}
    m_1,\dots,m_d=1 \\
    \mu_1,\dots,\mu_\delta\neq 1 \\
\end{array}%
\right.$$ is coherent. If this is the case, then it has a solution
$(s_1,\dots,s_n)$ in $\Omega$. Thus some group in $\mathcal{S}$,
namely the subgroup of $\Omega$ generated by $(s_1,\dots,s_n)$,
embeds in $\Omega$. This proves that the skeleton of $\Omega$ is
dense.\nolinebreak
\end{proof}

Propositions \ref{prop word_not_dense} and \ref{prop skeleton
dense} together prove that every existentially closed group
contains a finitely generated subgroup with unsolvable word
problem. On the other hand, Macintyre \cite{MR0317928} (see also
\cite[Chap.~IV, Theorem 8.5]{MR0577064}) has proved that there
exist two existentially closed groups $\Omega_1,\Omega_2$ such
that the intersection of the skeletons of the two is reduced to
the set of groups with solvable word problem. Moreover, Boone and
Higman \cite{MR0357625} (see also \cite[Chap.~IV, Theorem
7.4]{MR0577064}) have proved that every group with solvable word
problem embeds in a simple subgroup of a finitely presented group;
this easily implies (\cite{MR0414671} or \cite[Chap.~IV, Theorem
8.4]{MR0577064}) that a finitely generated group with solvable
word problem embeds in every existentially closed group.
We leave open the following question

\begin{question}Does every finitely generated group with solvable word
problem embed into an isolated group?
\end{question}

Note that the stronger well-known question whether every finitely
generated group with solvable word problem embeds in a finitely
presented simple group is open.

%%%%%%%%%%%%%%%%%%%%%%%%%%%%%%%%%%%%%%%%%%%%%%%%%%%%%%%%%%%%%%%%%%%%%%
%%%%%%%%%%%%%%%%%%%%%%%%%%%%%%%%%%%%%%%%%%%%%%%%%%%%%%%%%%%%%%%%%%%%%%
\section{Hereditary constructions for isolated groups}\label{sec hereditary}
%%%%%%%%%%%%%%%%%%%%%%%%%%%%%%%%%%%%%%%%%%%%%%%%%%%%%%%%%%%%%%%%%%%%%%
%%%%%%%%%%%%%%%%%%%%%%%%%%%%%%%%%%%%%%%%%%%%%%%%%%%%%%%%%%%%%%%%%%%%%%

The two characterizing properties for isolated points, finite
presentation and finite discriminability, are quite different in
nature.  The hereditary problem for finite presentation is
classical, and in most cases is well understood. Therefore our
results mostly deal with the hereditary problem for finite
discriminability.

%%%%%%%%%%%%%%%%%%%%%%%%%%%%%%%%%%%%%%%%%%%%%%%%%%%%%%%%%%%%%%%%%%%%%%
%%%%%%%%%%%%%%%%%%%%%%%%%%%%%%%%%%%%%%%%%%%%%%%%%%%%%%%%%%%%%%%%%%%%%%
%\subsection{Extensions of groups}\label{subsec extensions}
%%%%%%%%%%%%%%%%%%%%%%%%%%%%%%%%%%%%%%%%%%%%%%%%%%%%%%%%%%%%%%%%%%%%%%
%%%%%%%%%%%%%%%%%%%%%%%%%%%%%%%%%%%%%%%%%%%%%%%%%%%%%%%%%%%%%%%%%%%%%%

The analysis of extensions of isolated groups rests on the
analysis of their centre. As the centre of a finitely
discriminable group is itself finitely discriminable (see Lemma
\ref{lem fd centre} below), the analysis of the centre fits into
the more general problem of understanding finitely discriminable
abelian groups.

A group $G$ is finitely cogenerated if it has a finite subset $F$
having non-empty intersection with every non-trivial subgroup of
$G$.

The class of finitely cogenerated groups is clearly contained in
the class of finitely discriminable groups, but is much smaller;
for instance, it is closed under taking subgroups, and therefore a
finitely cogenerated group is necessarily torsion. However in
restriction to abelian groups, the two classes obviously coincide.

We denote by $C_{p^\infty}$ the $p$-primary Pr\"ufer group (also
called quasi-cyclic): this is the direct limit of cyclic groups of
order $p^n$ when $n\to\infty$; it can directly be constructed as
the quotient group $\Z[1/p]/\Z$.

In \cite{MR0179245} Yahya characterizes finitely cogenerated (i.e.
finitely discriminable) abelian groups.

\begin{lemma}\label{lem centerofisolated}
For an abelian group $G$, the following are equivalent
\begin{itemize}
    \item[(i)] $G$ is finitely discriminable.
    \item[(ii)] $G$ is artinian: every descending sequence of subgroups
    stabilizes.
    \item[(iii)] The three following conditions are satisfied
    \begin{enumerate}
        \item $G$ is a torsion group.
        \item Its $p$-torsion $\{x\in G\,|\;px=0\}$ is finite for all primes~$p$.
        \item $G$ has non-trivial $p$-torsion for only
        finitely many primes~$p$.
\end{enumerate}
    \item[(iv)] $G$ is a finite direct sum of finite cyclic groups and
    Pr\"ufer groups.\nolinebreak\qed
\end{itemize}
\end{lemma}

Lemma \ref{lem centerofisolated} is useful even when we focus on
finitely generated groups, in view of the following fact.

\begin{lemma}\label{lem fd centre}If $G$ is a finitely discriminable group, then its
centre is finitely discriminable.
\end{lemma}
\begin{proof}This immediately follows from the fact that every
subgroup of the centre of $G$ is normal in $G$.
\end{proof}

The converse of Lemma \ref{lem fd centre} is true in the case of
nilpotent groups, and more generally hypercentral groups. Recall
that, in a group $G$, the transfinite ascending central series
$(Z_\alpha)$ is defined as follows: $Z_1$ is the centre of $G$,
$Z_{\alpha+1}$ is the preimage in $G$ of the centre of
$G/Z_\alpha$, and $Z_\lambda=\bigcup_{\alpha<\lambda}Z_\alpha$ if
$\lambda$ is a limit ordinal. The group $G$ is
\textit{hypercentral} if $Z_\alpha=G$ for some $\alpha$.

%%%%%%%%%%%%%%%%%%%%%%%%%%%%%%%%%%%%%%%%%%%%%%%%%%%%%%%%%%%%%%%%%%%%%
\begin{corollary}\label{cor nilpotent}
 %%%%%%%%%%%%%%%%%%%%%%%%%%%%%%%%%%%%%%%%%%%%%%%%%%%%%%%%%%%%%%%%%%%%%
Let $G$ be a hypercentral group. Then $G$ is finitely
discriminable if and only if its centre $Z(G)$ is so.
%In particular an isolated nilpotent group is finite.
\end{corollary}
%%%%%%%%%%%%%%%%%%%%%%%%%%%%%%%%%%%%%%%%%%%%%%%%%%%%%%%%%%%%%%%%%%%%%
\begin{proof}
As noticed before, the ``if" part is straightforward. The converse
implication follows immediately from the known fact that if $G$ is
hypercentral, then any normal subgroup of $G$ intersects the
centre $Z(G)$ non-trivially. Let us recall the argument. Suppose
that a normal subgroup $N\neq 1$ has trivial intersection with the
centre. Let $\alpha$ be the smallest ordinal such that $Z_\alpha$
contains a non-trivial element $x$ of $N$. Clearly, $\alpha$ is a
successor. Let $M$ be the normal subgroup generated by $x$. Then
$M\subset N\cap Z_\alpha$. On the other hand, $M\cap
Z_{\alpha-1}=1$, and since $M\subset Z_\alpha$, by definition of
$Z_\alpha$ we have $[G,Z_\alpha]\subset Z_{\alpha-1}$. Hence
$[G,M]\subset Z_{\alpha-1}\cap M=1$, so that $M$ is central, a
contradiction.
\end{proof}

We now study hereditary properties of finitely discriminable
groups. It is convenient to extend some of the definitions above.
If $G$ is a group, we define a $G$-group as a group $H$ endowed
with an action of $G$ (when $H$ is abelian it is usually called a
$G$-module). For instance, normal subgroups and quotients of $G$
are naturally $G$-groups. We call a $G$-group finitely
discriminable (or $G$-finitely discriminable) if $\{1\}$ is
isolated among normal $G$-subgroups of $H$. (Proposition \ref{thm
finimin} has an obvious analog in this context.)

%%%%%%%%%%%%%%%%%%%%%%%%%%%%%%%%%%%%%%%%%%%%%%%%%%%%%%%%%%%%%%%%%%%%%
\begin{theorem}\label{thm stabextension}
%%%%%%%%%%%%%%%%%%%%%%%%%%%%%%%%%%%%%%%%%%%%%%%%%%%%%%%%%%%%%%%%%%%%%
Consider an extension of groups
$$1\loto K\loto G\loto Q\loto 1.$$
Denote by $W$ the kernel of the natural homomorphism
$Q\to\Out(K)$.

\begin{itemize}
\item[(1)] $G$ is finitely discriminable if and only if $K$ and
$C_G(K)$ are both $G$-finitely discriminable.

\item[(2)] Suppose that $K$ and $W$ are both $G$-finitely
discriminable. Suppose moreover that

\begin{itemize}\item[(*)] $Z(K)$ contains no infinite simple $G$-submodule
$G$-isomorphic to a normal subgroup of $Q$ contained in
$W$.\end{itemize}

\noindent Then $G$ is finitely discriminable.
\end{itemize}
\end{theorem}
\begin{remark}1) Assumption (*) is satisfied when $Z(K)$ or $W$
does not contain any infinite-dimensional vector space over a
prime field ($\mathbf{F}_p$ or $\mathbf{Q}$). In particular, if
$Z(K)$ is artinian it is satisfied.

2) In Assumption (*), we can replace ``infinite" by ``with
infinite endomorphism ring", which is a more natural hypothesis in
view of the subsequent proof.

3) Assumption (*) cannot be dropped: see the example in Remark
\ref{rem not_finite_index}.
\end{remark}

\begin{proof}[Proof of Theorem \ref{thm stabextension}(1)]
The conditions are clearly necessary. Conversely suppose that they
are satisfied. Let $N_i$ be a net of normal subgroups of $G$
tending to~1. Then $N_i\cap K\to 1$; this implies that eventually
$N_i\cap K=1$. Thus eventually $N_i\subset C_G(K)$. But similarly
eventually $N_i\cap C_G(K)=1$. Accordingly eventually $N_i=1$.
\end{proof}

%%%%%%%%%%%%%%%%%%%%%%%%%%%%%%%%%%%%%%%%%%%%%%%%%%%%%%%%%%%%%%%%%%%%%
Before proving (2), we need some preliminary results. Consider an
extension of groups:
$$1\loto K\loto G\stackrel{\pi}\loto Q\loto 1.$$
For any normal subgroup $H \triangleleft Q$, denote by
$\mathcal{I}(H)$ the set of normal subgroups of $G$ that are sent
isomorphically onto $H$ by the projection~$\pi$.

%%%%%%%%%%%%%%%%%%%%%%%%%%%%%%%%%%%%%%%%%%%%%%%%%%%%%%%%%%%%%%%%%%%%
\begin{lemma}\label{lem rappel}
%%%%%%%%%%%%%%%%%%%%%%%%%%%%%%%%%%%%%%%%%%%%%%%%%%%%%%%%%%%%%%%%%%%%
Let $M$ be a normal subgroup in $G$ such that $M \cap K$ is
trivial and denote by $\pi(M)$ its image in $Q$. Then:
\begin{enumerate}
    \item[(1)] Any group in $\mathcal{I}(\pi(M))$ is isomorphic as a
    $G$-group to $\pi(M)$ via $\pi$ and in particular is naturally a
    $Q$-group.
    \item[(2)] If $\pi(M)$ is minimal in $Q$ then all groups
    in $\mathcal{I}(\pi(M))$ are minimal in $G$.
    \item[(3)] The set
    $\mathcal{I}(\pi(M))$  is in one-to-one correspondence with the
    set of  $G$-equivariant group homomorphisms $\textnormal{Hom}_G
    (\pi(M),Z(K))$.\\
\end{enumerate}
\end{lemma}
%%%%%%%%%%%%%%%%%%%%%%%%%%%%%%%%%%%%%%%%%%%%%%%%%%%%%%%%%%%%%%%%%%%%
\begin{proof}

Points $(1)$ and $(2)$ are easy exercises; we concentrate on the
proof of the third point.

For any element $H\in \mathcal{I}(\pi(M))$, denote by $\sigma_{H}$
the inverse homomorphism to $\pi: H \rightarrow \pi(M)$. There is
an obvious one-to-one correspondence between the set
$\mathcal{I}(\pi(M))$, and the set of homomorphisms $\{ \sigma_H \
\vert \ H \in \mathcal{I}(\pi(M)) \}$ given by the map $H
\longmapsto \sigma_H$ and its inverse
$\sigma_H\longmapsto\text{Image}(\sigma_H)$. We use these maps to
identify these two sets.

We claim that there is a faithful transitive action of the abelian
group \allowbreak$\text{Hom}_G (\pi(M), Z(K))$ on the latter set.
For $\phi \in  \text{Hom}_G (\pi(M), Z(K))$ and $\sigma_H \in
\mathcal{I}(\pi(M))$ we define
\[
\begin{array}{rcl}
\phi\cdot\sigma_H:\;\pi(M) & \longrightarrow & G \\
x & \longmapsto & \sigma_{H}(x) \phi_(x).
\end{array}
\]

We first prove that $\phi\cdot\sigma_H $ is indeed a group
homomorphism. If $x,y  \in \pi(M)$, then
\begin{eqnarray*}
(\phi\cdot\sigma_H)(xy) & = & \sigma_H(xy) \phi(xy) \\
            & = & \sigma_H(x) \sigma_H(y) \phi(x) \phi(y) \\
            & = & \sigma_H (x) \phi(x) [\phi(x)^{-1},\sigma_H(y)] \sigma_H (y)
            \phi(y).
\end{eqnarray*}

Since $H$ and $Z(K)$ are both normal subgroups of $G$, we have
$[\phi(x) ^{-1},\sigma_H (y)]\in [Z(K), H]\subset Z(K)\cap
H=\{1\}$, so that
$$(\phi\cdot\sigma_H)(xy)=(\phi\cdot\sigma_H)(x)(\phi\cdot\sigma_H)(y).$$

As $\sigma_H$ and $\phi$ are $G$-equivariant homomorphisms, it is
immediate that $\phi\cdot\sigma_H$ is $G$-equivariant; it is also
immediate that it defines an action of $\text{Hom}_G (\pi(M),
Z(K))$ on $\mathcal{I}(\pi(M))$.

To see that the action is transitive, we consider an element
$\sigma_H$  in $\mathcal{I}(\pi(M))$ and we define the
``transition map" from $\sigma_M$ to $\sigma_H$:
\[
\begin{array}{rcl}
\phi_{H}:  \pi(M) & \rightarrow & G \\
x &  \mapsto & \sigma_M(x)^{-1} \sigma_{H}(x).
\end{array}
\]

We claim that the map \(\phi_{H}\) is a $G$-equivariant
homomorphism and has values in \( Z(K)\).

Indeed, by construction \(\pi \circ \phi_{H}(x)=1\) for all \(x\)
in \(\pi(M)\), thus the image of \(\phi_{H}\) is contained in
\(K\). More precisely \(\phi_{H}(\pi(M)) \subset HM\cap K\). Since
$H$ and $M$ are normal subgroups and have trivial intersection
with $K$, they are contained in the centralizer of $K$ in $G$. It
follows that $\phi_{H}(\pi(M))$ is contained in the centre of $K$.
The fact that $\phi_H$ is a $G$-equivariant homomorphism follows
now from direct computation.

Finally the stabilizer of $\sigma_M$ is trivial, as $\phi\cdot\sigma_M
 = \sigma_M$ if and only if for all $x \in \pi(M)$
$\sigma_M (x) \phi(x) = \sigma_M(x)$, that is, $\phi$ is
identically 1. This ends the proof of Lemma \ref{lem rappel}.
\end{proof}

\begin{proof}[Proof of Theorem \ref{thm stabextension}(2)]
Denote by $L_1,\dots,L_d$ the minimal normal subgroups of $G$
contained in $K$. Let $N$ be a normal subgroup of $G$. If $N\cap
K\neq 1$, then $N$ contains some $L_i$.

Let us assume now that $N\cap K=1$. Then $N$ is contained in the
centralizer $C_G(K)$. On the other hand, one can check that the
image of $C_G(K)$ in $Q$ is $W$. Denote by $Q_1,\dots,Q_n$ the
minimal normal subgroups of $Q$ contained in $W$. Then the image
of $N$ in $Q$ contains some $Q_i$. It follows from Lemma \ref{lem
rappel}(2) that $\pi^{-1}(Q_i)\cap N$ is a minimal normal subgroup
of $G$ (belonging to $\mathcal{I}(Q_i)$). Thus it remains to prove
that $\mathcal{I}(Q_i)$ is finite for every $i$. By Lemma \ref{lem
rappel}(3), if non-empty, this set is in one-to-one correspondence
with \allowbreak$\text{Hom}_G(Q_i,Z(K))$; let us show that this
set is finite. We discuss the possible cases.

\begin{itemize}
    \item If $Q_i$ is non-abelian, then, since it is
characteristically simple, it is perfect, and therefore
$\text{Hom}_G(Q_i,Z(K))\subset \text{Hom}(Q_i,Z(K))=\{1\}$.
    \item If $Q_i$ is infinite abelian, then by the assumption (*), we
know that the centre $Z(K)$ does not contain any $G$-submodule
isomorphic to $Q_i$, and therefore $\text{Hom}_G(Q_i,Z(K))=\{1\}$.
    \item Suppose that $Q_i$ is finite abelian. Let $V$ be the sum of all
$G$-submodules of $Z(K)$ isomorphic to $Q_i$. Since $Q_i$ is a
simple $G$-module, one can check that $V$, as a $G$-module, is a
direct sum $\bigoplus_{j\in J}V_j$ of submodules $V_j$ isomorphic
as $G$-modules to $Q_i$. Since $Z(K)$ is $G$-finitely
discriminable, the index set $J$ must necessarily be finite.
Therefore $\text{Hom}_G(Q_i,Z(K))=\text{Hom}_G(Q_i,\bigoplus_{j\in
I}V_j)\simeq \prod_{j\in J}\text{End}_G(Q_i)$, which is
finite.\qedhere
\end{itemize}
\end{proof}

\begin{corollary}\label{cor extensions_fd}
The classes of finitely discriminable groups and of isolated
groups are closed under extensions.
\end{corollary}
\begin{proof}Since $K$ is finitely discriminable, it is $G$-finitely
discriminable, and since $Q$ is finitely discriminable, its normal
subgroup $W$ must be $Q$-finitely discriminable. We can therefore
apply Theorem \ref{thm stabextension}(2), noting that since $K$ is
finitely discriminable, its centre is artinian (Lemmas \ref{lem
centerofisolated} and \ref{lem fd centre}). It follows that the
assumption (*) is satisfied: indeed, if a subgroup of $Z(K)$ is a
simple $G$-submodule, then it must be contained in the $p$-torsion
of $Z(K)$ for some prime $p$ and therefore is finite.

The second assertion is obtained by combining this with the fact
that the class of finitely presented groups is closed under
extensions.
\end{proof}

\begin{corollary}\label{cor ext_out}
Consider an extension of groups
$$1\loto K\loto G\loto Q\loto 1.$$
Suppose that $K$ is $G$-finitely discriminable, and that the
natural homomorphism $Q\to\Out(K)$ is injective. Then $G$ is
finitely discriminable.
\end{corollary}
\begin{proof}Note that since the kernel $W$ of the
natural homomorphism $Q\to\Out(K)$ is trivial, the assumption (*)
of Theorem \ref{thm stabextension} is trivially
satisfied.\end{proof}

\begin{corollary}\label{cor fd_overgroups_fi}
To be finitely discriminable and to be isolated are properties
inherited by overgroups of finite index.
\end{corollary}
\begin{proof}
Let $G$ be a group, and $H$ a finitely discriminable subgroup of
finite index. Let $N$ be a subgroup of finite index of $H$ which
is normal in $G$. Then $N$ is $H$-finitely discriminable and
therefore $G$-finitely discriminable. Since $G/N$ is finite, it is
clear that all assumptions of Theorem \ref{thm stabextension}(2),
are satisfied, so that $G$ is finitely discriminable.
\end{proof}

\begin{remark}In contrast, the finite discriminable property does not pass
to subgroups of finite index, as the following example shows.
Consider the wreath product $\Gamma=\Z \wr \Z$. In
\cite{MR0110750} P.~Hall has constructed that there exists a
simple faithful $\Z\Gamma$-module $V$ whose underlying abelian
group is an infinite dimensional $\mathbf{Q}$-vector space. Then
there is an obvious action of $\Gamma\times \Z/2\Z$ on the direct
sum $V \oplus V$, where the cyclic group permutes the two copies.
Consider then the semi-direct product $G = (V\oplus V) \rtimes
(\Gamma\times \Z/2\Z)$. Its subgroup of index two $(V \oplus V)
\rtimes \Gamma$ is not finitely discriminable as for each rational
$r$ we have a different minimal subgroup $V_r = \{ (v,rv) \in
V\oplus V \ | \ v \in V\}$. Nevertheless, $G$ is finitely
discriminable since every non-trivial normal subgroup of $G$
contains one of its two minimal normal subgroups $V_1$ and
$V_{-1}$.

%The only minimal normal subgroup of $G$ is $V\oplus V$ and is
%contained in every non-trivial normal subgroup of $G$ because the
%action is faithful; therefore $G$ is finitely discriminable. But the

Remark that $G$ is not finitely presented. Indeed, $\Z\wr\Z$ is
not a quotient of any finitely presented solvable group
\cite{MR0120269,MR0591649}. However, we conjecture that there
exists an example of an isolated group having a non-isolated
finite index subgroup. \label{rem not_finite_index}\end{remark}

%%%%%%%%%%%%%%%%%%%%%%%%%%%%%%%%%%%%%%%%%%%%%%%%%%%%%%%%%%%%%%%%%%%%%%
%%%%%%%%%%%%%%%%%%%%%%%%%%%%%%%%%%%%%%%%%%%%%%%%%%%%%%%%%%%%%%%%%%%%%%
\section{Examples of isolated groups}\label{sect exaisolated}
%%%%%%%%%%%%%%%%%%%%%%%%%%%%%%%%%%%%%%%%%%%%%%%%%%%%%%%%%%%%%%%%%%%%%%
%%%%%%%%%%%%%%%%%%%%%%%%%%%%%%%%%%%%%%%%%%%%%%%%%%%%%%%%%%%%%%%%%%%%%%
%%%%%%%%%%%%%%%%%%%%%%%%%%%%%%%%%%%%%%%%%%%%%%%%%%%%%%%%%%%%%%%%%%%%%%
\subsection{Elementary class}\label{subs:elementary_class}
%%%%%%%%%%%%%%%%%%%%%%%%%%%%%%%%%%%%%%%%%%%%%%%%%%%%%%%%%%%%%%%%%%%%%%
Obvious examples of isolated groups are finitely presented simple
groups. According to Corollary \ref{cor extensions_fd}, to get new
examples of isolated groups it suffices to consider the class of
groups that can be obtained from these by successively taking
extensions of groups. The class of groups we get is fairly well
understood, for by definition any such group has a composition
series of finite length, and by a theorem of Wielandt
\cite{MR1545814} this is precisely the class of finitely presented
groups that contain finitely many subnormal subgroups. More
generally we have

%%%%%%%%%%%%%%%%%%%%%%%%%%%%%%%%%%%%%%%%%%%%%%%%%%%%%%%%%%%%%%%%%%%%%%
\begin{proposition}
%%%%%%%%%%%%%%%%%%%%%%%%%%%%%%%%%%%%%%%%%%%%%%%%%%%%%%%%%%%%%%%%%%%%%%
Any group with finitely many normal subgroups is finitely
discriminable. In particular, finitely presented groups with
finitely many normal subgroups are isolated.\qed
\end{proposition}
%%%%%%%%%%%%%%%%%%%%%%%%%%%%%%%%%%%%%%%%%%%%%%%%%%%%%%%%%%%%%%%%%%%%%%

\begin{remark}
A finitely presented group $G$ has finitely many normal subgroups
if and only if every quotient of $G$ is isolated. The condition is
clearly necessary. Conversely if $G$ has infinitely many normal
subgroups, then, by compactness of $\mathcal{G}(G)$, there exists
an accumulation point, and hence $G$ has a non-isolated quotient.

Note, in contrast, that the Pr\"ufer group $C_{p^\infty}$ has all
its quotients finitely discriminable but has infinitely many
normal subgroups.
\end{remark}

\subsection{Quotients of isolated groups}\label{subs:quotient_isolated}
The inclusion above is strict in a strong sense:

\begin{theorem}\label{thm quot_isolated}
Every finitely generated group is quotient of an isolated group.
\end{theorem}

The theorem easily follows from the following lemma, which is
probably known but for which we found no reference.

\begin{lemma}\label{lem isolated out_free}
There exists an isolated group $K$ such that $\Out(K)$ contains a
non-abelian free group.
\end{lemma}

\begin{proof}[Proof of Theorem \ref{thm quot_isolated}]
Clearly it suffices to deal with a free group $F_n$. Consider $K$
as in Lemma \ref{lem isolated out_free}. Then $\Out(K)$ contains a
free group of rank $n$. Lift it to $\Aut(K)$ and consider the
semidirect product $G=K\rtimes F_n$ given by this action. By
Corollary \ref{cor ext_out}, $G$ is isolated.
\end{proof}

We prove later Lemma \ref{lem isolated out_free} (see Proposition
\ref{prop varabels}).

%%%%%%%%%%%%%%%%%%%%%%%%%%%%%%%%%%%%%%%%%%%%%%%%%%%%%%%%%%%%%%%%%
\subsection{Houghton groups}\label{subs:houghton}
%%%%%%%%%%%%%%%%%%%%%%%%%%%%%%%%%%%%%%%%%%%%%%%%%%%%%%%%%%%%%%%%%

These groups were first introduced by Houghton \cite{MR521478} in
his study of the relationship between ends and the cohomology of a
group. They were then studied by K. Brown in connection with the
so-called $FP$ cohomological properties \cite{MR885095}.

Fix an integer $n \geq 1$, let $\N$ denote the set of positive
integers and let $S = \N \times \{ 1,\dots, n\}$ denote a disjoint
union of $n$ copies of $\N$. Let $H_n$ be the subgroup of all
permutations $g$ of $S$ such that on each copy of $\N$, $g$ is
eventually a translation. More precisely, $g\in H_n$ if there is
an $n$-tuple $(m_1, \dots, m_n) \in \Z^n$ such that for each $i
\in \{ 1,\dots n\}$ one has $g(x,i) = (x + m_i,i)$ for all
sufficiently large $ x \in \N$. The map $g \mapsto (m_1, \dots,
m_n)$ is a homomorphism $\phi: H_n \rightarrow \Z^n$ whose image
is the subgroup $\{(m_1, \dots, m_n) \in \Z^n \ \vert \ \sum m_i
=0 \}$, of rank $n-1$. The kernel of $\phi$ is the infinite
symmetric group, consisting of all permutations of $S$ with finite
support. It coincides with the commutator subgroup of $H_n$ for $n
\geq 3$, while for $n = 1$ and $n = 2$, the commutator subgroup is
the infinite alternating group $\textnormal{Alt}(S)$. In all cases
the second commutator is the infinite alternating group, which is
a locally finite, infinite simple group.

\begin{proposition}\label{prop houghton_fd}
For every $n\ge 1$, the group $H_n$ is finitely discriminable.
\end{proposition}
\begin{proof}
There is an extension
$$1\loto\textnormal{Alt}(S)\loto H_n\loto
H_n/\textnormal{Alt}(S)\loto 1.$$ Since $\textnormal{Alt}(S)$ has
trivial centralizer in the full group of permutations of $S$, the
assumption of Corollary \ref{cor ext_out} is satisfied. Actually
$\textnormal{Alt}(S)$ is the unique minimal normal subgroup of
$H_n$ and is contained in all other non-trivial normal subgroups.
\end{proof}

The group $H_1$ is the infinite symmetric group and hence is not
finitely generated. The group $H_2$ is finitely generated, but not
finitely presented; indeed, it is a classical example of a
non-residually finite group that is a limit of finite groups
\cite{MR1411234,MR1458419}. For $n\ge 3$, it is a result of
K.~Brown \cite{MR885095} that $H_n$ is finitely presented.
Therefore by Proposition \ref{prop houghton_fd} it is isolated.

%%%%%%%%%%%%%%%%%%%%%%%%%%%%%%%%%%%%%%%%%%%%%%%%%%%%%%%%%%%%%%%%%
\subsection{Abels groups}\label{subs:abels}
%%%%%%%%%%%%%%%%%%%%%%%%%%%%%%%%%%%%%%%%%%%%%%%%%%%%%%%%%%%%%%%%%

Fix an integer $n \geq 2$ and a prime $p$. Denote by $A_n \subset
\textnormal{GL}_n(\Z[1/p])$ the subgroup of upper triangular
matrices $a$ such that $a_{11} = a_{nn} = 1$ and such that the
other diagonal coefficients are positive. For example
\[
A_2 = \left( \begin{matrix}
1 & \Z[1/p] \\
0 & 1
      \end{matrix}\right);\quad
A_3 = \left( \begin{matrix}
1 & \Z[1/p] & \Z[1/p] \\
0 & p^{\Z} & \Z[1/p] \\
0 & 0& 1
             \end{matrix} \right).
\]

The group $A_3$, introduced by Hall in \cite{MR0124406} is not
finitely presented \cite{MR0577069}. Abels \cite{MR564423}
introduced $A_4$ and showed that it is finitely presented, a
result that was subsequently extended for $n \geq 4$, see
\cite{MR885095} for a proof and a discussion of homological
properties of these groups.

The centre of $A_n$ consists of unipotent matrices with a single
possibly non-trivial element in the upper right corner. It is
clearly isomorphic to $\Z[1/p]$. The conjugation by the diagonal
matrix $\text{Diag}(p,1,\dots,1)$ provides an automorphism of
$A_n$ which induces the multiplication by $p$ on the centre.
Consider the canonical copy $Z$ of $\Z$ contained in the centre
through its identification with $\Z[1/p]$. Then the quotient
$A_n/Z$ is non-Hopfian (see \cite{MR564423} for details). As
noticed before, we have

%%%%%%%%%%%%%%%%%%%%%%%%%%%%%%%%%%%%%%%%%%%%%%%%%%%%%%%%%%%%%%%%%
\begin{proposition}\label{prop abelsfinpres}
%%%%%%%%%%%%%%%%%%%%%%%%%%%%%%%%%%%%%%%%%%%%%%%%%%%%%%%%%%%%%%%%%
The groups $A_n$ and $A_n/Z$ are finitely presented for all $n
\geq 4$.\qed
\end{proposition}
%%%%%%%%%%%%%%%%%%%%%%%%%%%%%%%%%%%%%%%%%%%%%%%%%%%%%%%%%%%%%%%%%

We now turn to finite discriminability. The group $A_n$ itself is
certainly not finitely discriminable because its centre $\Z[1/p]$
is not an artinian abelian group (or because it is residually
finite). In contrast, the centre of $A_n/Z$ is a Pr\"ufer group
$C_{p^\infty}$.

%%%%%%%%%%%%%%%%%%%%%%%%%%%%%%%%%%%%%%%%%%%%%%%%%%%%%%%%%%%%%%%%%
\begin{proposition}\label{prop abelsfindiscr}
%%%%%%%%%%%%%%%%%%%%%%%%%%%%%%%%%%%%%%%%%%%%%%%%%%%%%%%%%%%%%%%%%
The groups $A_n/Z$ are finitely discriminable for $n \geq 2$. In
particular, for $n \geq 4$ these groups are infinite, solvable
(3-solvable when $n=4$), and isolated.
\end{proposition}
%%%%%%%%%%%%%%%%%%%%%%%%%%%%%%%%%%%%%%%%%%%%%%%%%%%%%%%%%%%%%%%%%

Before giving the proof of the proposition, let us give some of
its consequences. First notice that this allows to obtain a kind
of converse of Lemma \ref{lem fd centre}.

%%%%%%%%%%%%%%%%%%%%%%%%%%%%%%%%%%%%%%%%%%%%%%%%%%%%%%%%%%%%%%%%%
\begin{corollary}
%%%%%%%%%%%%%%%%%%%%%%%%%%%%%%%%%%%%%%%%%%%%%%%%%%%%%%%%%%%%%%%%%
Every finitely discriminable abelian group is isomorphic to the
centre of an isolated group.
\end{corollary}
%%%%%%%%%%%%%%%%%%%%%%%%%%%%%%%%%%%%%%%%%%%%%%%%%%%%%%%%%%%%%%%%%
\begin{proof}In this proof, let us denote the Abels group
by $A_{n,p}$ to make explicit the dependance on $p$. Let $G$ be a
finitely discriminable abelian group. By Lemma \ref{lem
centerofisolated}, $G$ is isomorphic to
$F\times\prod_{i=1}^nC_{{p_i}^\infty}$, where $F$ is a finite
abelian group, and $p_1,\dots,p_n$ are prime. Then $G$ is
isomorphic to the centre of
$F\times\prod_{i=1}^nA_{4,p_i}/Z(A_{4,p_i})$, which is an isolated
group by Proposition \ref{prop abelsfindiscr} and Corollary
\ref{cor extensions_fd}.
\end{proof}

%%%%%%%%%%%%%%%%%%%%%%%%%%%%%%%%%%%%%%%%%%%%%%%%%%%%%%%%%%%%%%%%%
\begin{corollary}\label{cor hopfdensty}
%%%%%%%%%%%%%%%%%%%%%%%%%%%%%%%%%%%%%%%%%%%%%%%%%%%%%%%%%%%%%%%%%
The Hopfian property is not dense in the space of finitely
generated groups.\qed
\end{corollary}
%%%%%%%%%%%%%%%%%%%%%%%%%%%%%%%%%%%%%%%%%%%%%%%%%%%%%%%%%%%%%%%%%

The non-Hopfian property is clearly not dense since finite groups
are isolated Hopfian. The Hopfian Property is not open
\cite{ABLRSV,Stalder}: the residually finite (metabelian) groups
$\Z\wr\Z$ and $\Z[1/6]\rtimes_{3/2}\Z$ are approximable by
non-Hopfian Baumslag-Solitar groups. On the other hand, let us
allow ourselves a little digression:

\begin{proposition}
The Hopfian Property is not closed in the space of finitely
generated groups. More precisely, there exists a finitely
generated solvable group that is approximable by finite groups,
but is non-Hopfian.
\end{proposition}
\begin{proof}
For $n\ge 3$, consider the group $B_n$ defined in the same way as
$A_n$, but over the ring $\mathbf{F}_p[t,t^{-1}]$ rather than
$\Z[1/p]$. It is easily checked to be finitely generated, and its
centre can be identified with $\mathbf{F}_p[t,t^{-1}]$, mapping on
the upper right entry. Similarly to the case of $A_n$, the group
$B_n/\mathbf{F}_p[t]$ is non-Hopfian. On the one hand, write
$\mathbf{F}_p[t]$ as the union of an increasing sequence of finite
additive subgroups $H_k$. Then $B_n/H_k$ converges to the
non-Hopfian group $B_n/\mathbf{F}_p[t]$ when $k\to\infty$. On the
other hand, as a finitely generated linear group, $B_n$ is
residually finite, and therefore so is $B_n/H_k$. So the
non-Hopfian finitely generated group $B_n/\mathbf{F}_p[t]$ is
approximable by finite groups.
\end{proof}

Before proving Proposition \ref{prop abelsfindiscr}, let us
describe a variation of Abels' group introduced in
\cite{CornulierP}. Consider integers $n_1,n_2,n_3,n_4$ satisfying
$n_1,n_4\ge 1$ and $n_2,n_3\ge 3$. Consider the group $H$ of upper
triangular matrices by blocks $(n_1,n_2,n_3,n_4)$ of the form
$$ \begin{pmatrix}
  I_{n_1} & (*)_{12} & (*)_{13} & (*)_{14} \\
  0 & (**)_{22} & (*)_{23} & (*)_{24} \\
  0 & 0 & (**)_{33} & (*)_{34} \\
  0 & 0 & 0 & I_{n_4} \\
\end{pmatrix},$$
where $(*)$ denote any matrices and $(**)_{ii}$ denote matrices in
$\textnormal{SL}_{n_i}$, $i=2,3$.

Set $\Gamma=H(\Z[1/p])$, and let $Z$ be the subgroup consisting of
matrices of the form $I+A$ where the only nonzero block of $A$ is
the block $(*)_{14}$ and has integer coefficients; note that $Z$
is a free abelian group of rank $n_1n_4$, and is central in
$\Gamma$.

It is proved in \cite{CornulierP} that $\Gamma$ and $\Gamma/Z$
have Kazhdan's Property (T), are finitely presentable and that the
group $\textnormal{GL}_{n_1}(\Z)$ embeds in $\Out(\Gamma/Z)$.

\begin{proposition}\label{prop varabels}
The group $\Gamma/Z$ is isolated, has Kazhdan's Property (T), and
its outer automorphism group contains a non-abelian free subgroup
if $n_1\ge 2$.
\end{proposition}

It only remains to prove that $\Gamma/Z$ is finitely
discriminable. We now prove it together with Proposition \ref{prop
abelsfindiscr}.

\begin{proof}[Proof of Propositions \ref{prop abelsfindiscr} and \ref{prop varabels}]

Let us make a more general construction.

Consider integers $m_1,m_2,m_3\ge 1$, and $n=m_1+m_2+m_3$.
Consider the subgroup $U$ of $\textnormal{GL}_n$ given by upper
unipotent by $(m_1,m_2,m_3)$-blocks matrices, that is, matrices of
the form:$$\begin{pmatrix}
  I_{m_1} & A_{12} & A_{13} \\
  0 & I_{m_2} & A_{23} \\
  0 & 0 & I_{m_3} \\
\end{pmatrix}.$$
Let $V$ denote the subgroup of $U$ consisting of matrices with
$A_{12}=0$ and $A_{23}=0$.

\begin{lemma}\label{lem centralisateur_abels_gen}
The centralizer $C$ of $U$ modulo $V$ in $\textnormal{GL}_n$ is
reduced to $U_s$, the group generated by $U$ and scalar matrices.
\end{lemma}
\begin{proof}
First compute the normalizer $N$ of $U$ in $\textnormal{GL}_n$.
Denote by $E$ the vector space of rank $n$. Since the fixed points
of $U$ is the subspace $E_1$ generated by the $m_1$ first
coefficients, it must be invariant under $N$. Since the fixed
point of $U$ on $E/E_1$ is the subspace $E_2$ generated by the
$m_1+m_2$ first coefficients, it must also be invariant under $N$.
We thus obtain that $N$ is the group of upper triangular matrices
under this decomposition by blocks.

Let us now show that $C=U_s$. Since $U\subset C$, it suffices to
show that $C\cap D=S$, where $D$ is the group of diagonal by
blocks matrices and $S$ is the group of scalar matrices.

If we take $A=\begin{pmatrix}
  I_{m_1} & A_{12} & A_{13} \\
  0 & I_{m_2} & A_{23} \\
  0 & 0 & I_{m_3} \\
\end{pmatrix}\in U$ and $D=\begin{pmatrix}
  D_1 & 0 & 0 \\
  0 & D_2 & 0 \\
  0 & 0 & D_3 \\
\end{pmatrix}\in C\cap D$, then $[D,A]=\begin{pmatrix}
  I_{m_1} & D_1A_{12}{D_2}^{-1}-A_{12} & (\dots) \\
  0 & I_{m_2} & D_2A_{23}{D_3}^{-1}-A_{23} \\
  0 & 0 & I_{m_3} \\
\end{pmatrix}$ must belong to $V$. Thus $D_1A_{12}=A_{12}D_2$ and
$D_2A_{23}=A_{23}D_3$ for all $A_{12},A_{23}$. This easily implies
that there exists a scalar $\lambda$ such that $D_i=\lambda
I_{m_i}$ for each $i=1,2,3$.
\end{proof}

Denote now by $G$ the group of upper triangular by
$(m_1,m_2,m_3)$-blocks matrices with $A_{11}=I_{m_1}$ and
$A_{33}=I_{m_3}$. Note that $V$ is central in $G$.

\begin{lemma}\label{lem abels_gen_fd}
Let $R$ be a commutative ring. Let $H$ be a subgroup of $G(R)$
containing $U(R)$, and let $Z$ be any subgroup of $V(R)$
satisfying the following assumption: for every $x\in V(R)-\{0\}$,
there exists $\alpha\in R$ such that $\alpha x\notin Z$. Then
$H/Z$ is finitely discriminable if and only if the abelian group
$V(R)/Z$ is finitely discriminable.
\end{lemma}
\begin{proof}
By Lemma \ref{lem fd centre}, the condition is necessary since
$V(R)/Z$ is central in $H$.

Conversely, we have an extension $$1\to U(R)/Z\to H/Z\to H/U(R)\to
1.$$ By Lemma \ref{lem centralisateur_abels_gen}, the centralizer
of $U(R)/Z$ in $H/Z$ is contained in $U(R)/Z$, and therefore the
natural homomorphism $H/U(R)\to\Out(U(R)/Z)$ is injective. Now
$U(R)/Z$ is nilpotent, and by Corollary \ref{cor nilpotent} it is
finitely discriminable if and only if its centre is so. Thus it
suffices to prove that the centre of $U(R)/Z$ is $V(R)/Z$.

Suppose that a matrix $A=\begin{pmatrix}
  I_{m_1} & A_{12} & A_{13} \\
  0 & I_{m_2} & A_{23} \\
  0 & 0 & I_{m_3}\end{pmatrix}$
 is central in $U(R)/Z$. By an immediate
  computation it must satisfy, for all $B_{12},B_{23}$, the
  property $A_{12}B_{23}-B_{12}A_{13}\in Z$. If $A\notin V(R)$, we
  can choose $B_{12},B_{23}$ so that $x=A_{12}B_{23}-B_{12}A_{13}\neq
  0$. Choose $\alpha$ as in the assumption of the lemma. Then
  $A_{12}(\alpha B_{23})-(\alpha B_{12})A_{13}=\alpha x\notin Z$,
  a contradiction. Thus the centre of $U(R)/Z$ is $V(R)/Z$.
\end{proof}

In the case of Abels' group $A_n$, we have $R=\Z[1/p]$,
$m_1=m_3=1$, $m_2=n-2$, and $V(R)/Z$ is isomorphic to the Pr\"ufer
group $C_{p^{\infty}}$. The group $H$ is given by matrices in
$G(\Z[1/p])$ whose diagonal block 22 is upper triangular with
powers of $p$ on the diagonal. The assumption of Lemma \ref{lem
abels_gen_fd} is always satisfied, with $\alpha=p^{-k}$ for
some~$k$. This proves Proposition \ref{prop abelsfindiscr}.

In the case of the variant of Abels' group $\Gamma$, we have
$R=\Z[1/p]$, $m_1=n_1$, $m_2=n_2+n_3$, $m_3=n_4$, and $V(R)/Z$ is
isomorphic to ${(C_{p^\infty})}^{n_1n_4}$. The group $H$ is given
by matrices in $G(\Z[1/p])$ whose diagonal block 22 is upper
triangular by blocks $(n_2,n_3)$ with the two diagonal sub-blocks
of determinant one. The assumption of Lemma \ref{lem abels_gen_fd}
is also always satisfied, again with $\alpha=p^{-k}$ for some~$k$.
This proves Proposition \ref{prop varabels}.
\end{proof}

%%%%%%%%%%%%%%%%%%%%%%%%%%%%%%%%%%%%%%%%%%%%%%%%%%%%%%%%%%%%%%%%%
\subsection{Thompson groups}\label{subs:thompson}
%%%%%%%%%%%%%%%%%%%%%%%%%%%%%%%%%%%%%%%%%%%%%%%%%%%%%%%%%%%%%%%%%

The results we mention on these groups can be found in
\cite{MR1426438}. Let $F$ be the set of piecewise linear
increasing homeomorphisms from the closed unit interval $\lbrack
0,1 \rbrack$ to itself that are differentiable except at finitely
many dyadic rational numbers and such that on intervals of
differentiability the derivatives are powers of $2$. This turns
out to be a finitely presented group. There is a natural morphism
$F\nolinebreak\to\nolinebreak\Z^2$, mapping $f$ to $(m,n)$, where
the slope of $f$ at 0 is $2^m$ and the slope at 1 is $2^n$. The kernel
$F_0$ is a simple (infinitely generated) group. The corresponding
extension satisfies the assumptions of Corollary \ref{cor
ext_out}: more precisely, $F_0$ is contained in every non-trivial
normal subgroup of $F$. In particular

\begin{proposition}
Thompson's group $F$ is isolated.\qed
\end{proposition}

\noindent Observe the analogy with Houghton's group $H_3$
mentioned above.

Consider now $S^1$ as the interval $\lbrack 0,1 \rbrack$ with the
endpoints identified. Thompson's group $T$ is defined as the group
of piecewise linear homeomorphisms from $S^1$ to itself that map
images of dyadic rational numbers to images of dyadic rational
numbers and that are differentiable except at finitely many images
of rational dyadic numbers and on the intervals of
differentiability the derivatives are power of $2$. The group $T$
is finitely presented and simple, and in particular is isolated.
We use the groups $F$ and $T$ in the next paragraph.

%%%%%%%%%%%%%%%%%%%%%%%%%%%%%%%%%%%%%%%%%%%%%%%%%%%%%%%%%%%%%%%%%
\subsection{Wreath products}\label{subs:wreath}
%%%%%%%%%%%%%%%%%%%%%%%%%%%%%%%%%%%%%%%%%%%%%%%%%%%%%%%%%%%%%%%%%

If $G$ and $W$ are two groups, the \textit{standard wreath
product} $W\wr G$ is the semi-direct product $W^{(G)} \rtimes G$,
where $W^{(G)}$ denotes the direct sum of copies of $W$ indexed by
$G$, on which $G$ acts via the action on the labels by left
multiplication.

A wreath product $W\wr G$ of two finitely presented groups is
finitely presented only in trivial cases, namely when $G$ is
finite or $W=1$ (see \cite{MR0120269}). Nevertheless if we
consider \emph{permutational wreath products}, then positive
results on finite presentation do exist. Let $G$ and $W$ be groups
and $X$ a $G$-set, which we suppose transitive to simplify, so
that we can write $X=G/H$. The \emph{permutational wreath product}
$W\wr_X G$ is the semi-direct product $W^{(X)} \rtimes G$, where
$W^{(X)}$ denotes the direct sum of $X$ copies of $W$, on which
$G$ acts via the natural action on the labels. The following
theorem is shown in
\cite{Cornulier}.%

%%%%%%%%%%%%%%%%%%%%%%%%%%%%%%%%%%%%%%%%%%%%%%%%%%%%%%%%%%%%%%%%%%%%%%
\begin{theorem}\label{thm wreathfinpres}
%%%%%%%%%%%%%%%%%%%%%%%%%%%%%%%%%%%%%%%%%%%%%%%%%%%%%%%%%%%%%%%%%%%%%%
Let $G$ and $W\neq 1$ be groups, and $X=G/H$ a transitive $G$-set.
The group $W \wr_X G$ is finitely presented if and only if
\begin{enumerate}
    \item[(i)] both $W$ and $G$ are finitely presented;
    \item[(ii)] $H$ is finitely generated;
\item[(iii)] the product action of $G$ on $X \times X$ has
finitely many orbits (equivalently, the double coset space
$H\backslash G/H$ is finite).\qed
\end{enumerate}
\end{theorem}
%%%%%%%%%%%%%%%%%%%%%%%%%%%%%%%%%%%%%%%%%%%%%%%%%%%%%%%%%%%%%%%%%%%%%%%

\begin{proposition}\label{prop wreath_nocentre}
Keep the notation as above. Suppose that $X=G/H$ is a faithful
transitive $G$-set, and that $W\neq 1$ is finitely discriminable
and has trivial centre. Then the wreath product $W \wr_X G$ is
finitely discriminable.
\end{proposition}
\begin{proof}
Consider the extension
$$1\loto W^{(X)}\loto W \wr_X G\loto G\loto
1.$$ Since $W\neq 1$, the natural morphism $G\to\Out(W^{(X)})$ is
injective. So, to apply Corollary \ref{cor ext_out}, it suffices
to show that $W^{(X)}$ is a finitely discriminable $G$-group. Let
$N$ be a normal subgroup of $G$ contained in $W^{(X)}$. For $x\in
X$, denote by $W_x$ the $x$-th copy of $W$ in $W^{(X)}$. If $N\cap
W_x=1$, for some $x$, then by transitivity of the $G$-action on
$X$, we have $N\cap W_x=1$ for all $x$ and therefore $N$
centralizes all $W_x$. Since $W$ has trivial centre, this implies
$N=1$. Therefore if $N\neq 1$, then $N\cap W_x$ contains a minimal
normal subgroup $M$ of $W$. It follows that $N$ contains
$M^{(X)}$, and thus $W^{(X)}$ is a finitely discriminable
$G$-group.\nolinebreak\end{proof}

To deal with the case when $W$ has non-trivial centre we need some
further assumptions.

\begin{proposition}\label{prop wreath_centre}
Keep the notation as above. Suppose that $X=G/H$ is a faithful
$G$-set, and that $W\neq 1$ is finitely discriminable. Suppose
moreover that the two following conditions are
satisfied.\begin{itemize}
    \item[(i)] For every $n\ge 0$ and all distinct $y,x_1,\dots,x_n\in X$, there exists $g\in G$ such that $gy\neq y$
    and $gx_i=x_i$ for all $i=1,\dots,n$.
    \item[(ii)] The action of $G$ on $X^2$ has finitely many orbits.
\end{itemize}
Then the wreath product $W \wr_X G$ is finitely discriminable.
\end{proposition}
\begin{proof}
Arguing as in the proof of Proposition \ref{prop wreath_nocentre},
we are reduced to deal with a subgroup $N\neq 1$ contained in
$Z(W)^{(X)}$. Let $Z_p$ denote the $p$-torsion in $Z(W)$. Clearly,
for some $p$, ${Z_p}^{(X)}\cap N\neq 1$, and replacing $N$ by
$N\cap {Z_p}^{(X)}$ we can suppose $N$ contained in ${Z_p}^{(X)}$.

Consider a non-trivial element $w$ of $N$, and denote its support
as $\{y,x_1,\dots,x_n\}$. Using the assumption on the $G$ action
on $X$, there exists $g\in G$ fixing all $x_i$'s and mapping $y$
to some $y'\in X-\{y\}$. Then $[g,w]$ has support reduced to
$\{y,y'\}$. Consider a fixed finite family of elements $(u_i,v_i)$
in $X^2$ with an element in each $G$-orbit. There exists $i$ and
$h\in G$ mapping $(y,y')$ to $(u_i,v_i)$. Thus $h[g,w]h^{-1}$
belongs to $N$ and has support $\{u_i,v_i\}$. We obtain that if we
take elements with support some $\{u_i,v_i\}$ and with values in
elements of prime order in the centre of $W$, along with a finite discriminating subset in one copy of $W$, we obtain a finite
discriminating subset of $W\wr_X G$.
\end{proof}

Let us now give examples of $(G,X)$ satisfying the conditions of
both Theorem \ref{thm wreathfinpres} and Proposition \ref{prop
wreath_centre} (observe that the choice of the non-trivial
isolated group $W$ plays no role there). Trivial examples are
those when $X$ is finite.

\begin{itemize}
    \item {\it Houghton groups}. For $n\ge 1$, the group $H_n$ described
    in \S\ref{subs:houghton} acts on $X=\N\times\{1,\dots,n\}$. The action contains the
    groups of finitely supported permutations and therefore
    Assumption (iii) of Theorem \ref{thm wreathfinpres} and
    Assumptions (i) and (ii) of Proposition \ref{prop
    wreath_centre} are satisfied. The stabilizer of a point is also
    isomorphic to $H_n$. Accordingly, for $n\ge 3$, the group $H_n$ is
    finitely presented, the stabilizers are finitely generated,
    so that all assumptions are fulfilled.

    \item {\it Thompson groups}. Thompson's group $F$ (resp. $T$) acts
    on $X=\Z[1/2]\cap ]0,1[$ (resp. $X=\Z[1/2]/\Z$). This action
    is transitive on ordered (resp. cyclically ordered) $n$-tuples
    for all $n$. The stabilizer of a point is isomorphic to
    $F\times F$ (resp. $F$) and is therefore finitely generated.
    Thus all the assumptions of Theorem \ref{thm wreathfinpres} and
    Proposition \ref{prop wreath_centre} are satisfied.
\end{itemize}

%%%%%%%%%%%%%%%%%%%%%%%%%%%%%%%%%%%%%%%%%%%%%%%%%%%%%%%%%%%%%%%%%
\subsection{Grigorchuk's finitely presented amenable
group}\label{subs:grigorchuk_fp}
%%%%%%%%%%%%%%%%%%%%%%%%%%%%%%%%%%%%%%%%%%%%%%%%%%%%%%%%%%%%%%%%%

Let us consider the group $\widetilde{\Gamma}$ constructed by
Grigorchuk in \cite{MR1616436} and given by the presentation
$$\langle a,b,c,d,t \vert \, a^2=b^2=bcd=(ad)^4=(adacac)^4=1,$$
$$a^t=aca,b^t=d,c^t=b,d^t=c\rangle.$$

This group was provided as an example of a finitely presented
amenable group that is not elementary amenable \cite{MR1616436}.
The group $\widetilde{\Gamma}$ is an ascending HNN-extension over
the first Grigorchuk group $\Gamma=\langle a,b,c,d\rangle$,
introduced in \cite{MR0565099}, which, among other remarkable
properties, is the first known example of a group with
intermediate growth \cite{MR0764305}.

%%%%%%%%%%%%%%%%%%%%%%%%%%%%%%%%%%%%%%%%%%%%%%%%%%%%%%%%%%%%%%%%%
\begin{proposition}\label{prop isolgrig}
%%%%%%%%%%%%%%%%%%%%%%%%%%%%%%%%%%%%%%%%%%%%%%%%%%%%%%%%%%%%%%%%%
The group \(\widetilde{\Gamma}\) is isolated.
\end{proposition}
%%%%%%%%%%%%%%%%%%%%%%%%%%%%%%%%%%%%%%%%%%%%%%%%%%%%%%%%%%%%%%%%%
\begin{proof}
Sapir and Wise \cite{MR1868545} have proved that every proper
quotient of \(\widetilde{\Gamma}\) is metabelian. Therefore every
normal subgroup $N\neq 1$ contains the second derived subgroup of
\(\widetilde{\Gamma}\), which is not the trivial group since
\(\widetilde{\Gamma}\) itself is not metabelian.
\end{proof}

This result contradicts a conjecture of Stepin, appearing in
\cite[\S 1]{MR1616436} too, which states that the class of
elementary amenable finitely generated groups is dense in that of
amenable finitely generated groups.

\subsection{Deligne's central extension}\label{subs:deligne}

Deligne \cite{MR0507760} has shown that there exists a central
extension
$$1\loto Z\loto \widetilde{\Gamma}\loto\Gamma\loto 1,$$
where $\Gamma$ is a subgroup of finite index in
$\textnormal{PSp}_{2n}(\Z)$, the group $\widetilde{\Gamma}$ is its
preimage in the universal covering
$\widetilde{\textnormal{Sp}_{2n}(\R)}$ of
$\textnormal{PSp}_{2n}(\R)$, and $Z$ is infinite cyclic, that
satisfies the following remarkable property: every finite index
subgroup of $\widetilde{\Gamma}$ contains $Z$. By the
Kazhdan-Margulis Theorem \cite[Theorem 8.12]{MR0776417} every
non-trivial normal subgroup of $\Gamma$ has finite index. It
follows that every normal subgroup of $\widetilde{\Gamma}$ either
is contained in $Z$ or contains $Z$ (and has finite index in
$\widetilde{\Gamma}$). Thus $Z-\{1\}$ is a discriminating subset
for $\widetilde{\Gamma}$. This is infinite, but becomes finite
after taking a proper quotient of the centre. Thus we obtain:

\begin{proposition}
For every $n\ge 2$, the group $\widetilde{\Gamma}/nZ$ is
isolated.\qed
\end{proposition}

This shows that a (finite central) extension of an infinite
residually finite group can be isolated; moreover they are
lattices in the non-linear simple Lie group with finite centre
$\widetilde{\textnormal{Sp}_{2n}(\R)}/nZ$. Other similar examples
appear in \cite{MR0735524}, with $\Gamma$ a cocompact lattice in
$\textnormal{SO}(2,n)$ for $n\ge 3$. Erschler \cite{MR2029029}
provides examples where $\Gamma$ is the first Grigorchuk group;
these are examples of finitely generated, finitely discriminable
groups with intermediate growth.

\section{Further developements}

Let $K$ be a compact space. By induction on ordinals, define
$I_0(K)$ as the set of isolated points of $K$, and for $\alpha>0$,
define $I_{\alpha}(K)$ as the set of isolated points in
$K-\bigcup_{\beta<\alpha}I_\beta(K)$; call $I_\alpha(K)$ the set
of $\alpha$-isolated points of $K$. Let $\textnormal{Cond}(K)$
denote the condensation points, that is the complement of the
union $\bigcup_\alpha I_\alpha(K)$.

If $K$ is metrizable, then the sequence $(I_\alpha(K))$ breaks off
after a countable number of steps. In the case of the space of
finitely generated groups (in which case we simply write
$\mathbf{I}_\alpha$ and $\mathbf{Cond}$), what is this number?

It might be interesting to study $\mathbf{I}_\alpha$ for small
values of $\alpha$. For instance, it is easy to check that
$\Z^n\in \mathbf{I}_n$ for all $n$.

The study of $\mathbf{Cond}$ is also of interest. It is
characterized by: a finitely generated group $G$ is in
$\mathbf{Cond}$ if and only if every neighbourhood of $G$ is
uncountable. By the results of Champetier \cite{MR1760424}
$\mathbf{Cond}$ contains all non-elementary hyperbolic groups. It
can be showed that it also contains the wreath product $\Z\wr\Z$
and the first Grigorchuk group $\Gamma$.

A related variant of these definitions is the following: consider
a group $G$ (not necessarily finitely generated). We say that $G$
is in the class $\mathbf{II}_\alpha$ if $G\in
I_\alpha(\mathcal{G}(G))$. This holds for at most one $\alpha$;
otherwise say that $G$ belongs to the class $\mathbf{ICond}$. The
new $I$ stands for ``Inner" or ``Intrinsic". Note that, if we
restrict to finitely generated groups, we have
$\mathbf{ICond}\subset\mathbf{Cond}$, and every group in
$\mathbf{I}_\alpha$ belongs to $\mathbf{II}_\beta$ for some
$\beta\le\alpha$. Note also that the two classes coincide in
restriction to finitely presented groups, but for instance the
first Grigorchuk group $\Gamma$ belongs to $\mathbf{II}_1$.

\end{document}